\documentclass[leqno]{article}
\usepackage{amssymb}
\usepackage{amsmath, enumerate}
\usepackage{vmargin}

 \makeatletter
 \renewcommand\section{\@startsection {section}{1}{\z@}%
 {-3.5ex \@plus -1ex \@minus -.2ex}%
 {2.3ex \@plus.2ex}%
 {\center \normalfont\large\bfseries}}
 \makeatother

\setmarginsrb{3cm}{3cm}{3cm}{3cm}{1mm}{6mm}{0mm}{10mm}

\newtheorem{thm}{Theorem}[section]
\newtheorem{prop}[thm]{Proposition}
\newtheorem{cor}[thm]{Corollary}

\newtheorem{defi}[thm]{Definition}
\newtheorem{remark}[thm]{Remark}
\newtheorem{example}[thm]{Example}
\newtheorem{pb}[thm]{Problem}

\newenvironment{rk}{\begin{remark}\rm}{\end{remark}}
\newenvironment{definition}{\begin{defi}\rm}{\end{defi}}

\newcommand{\real}{{\mathbb R}}
\newcommand{\nat}{{\mathbb N}}

\newcommand{\T}{{\mathbb T}}

\newcommand{\A}{{\mathsf A}}
\newcommand{\B}{{\mathsf B}}
\newcommand{\D}{{\mathsf D}}
\newcommand{\F}{{\Phi}}
\newcommand{\M}{{\mathsf M}}

\renewcommand{\d}{\delta}
\newcommand{\De}{\Delta}
\newcommand{\e}{\varepsilon}

\renewcommand{\e}{\varepsilon}
\newcommand{\f}{\varphi}

\renewcommand{\o}{{\omega}}

\newcommand{\ot}{\otimes}

\renewcommand{\t}{\tau}
\newcommand{\8}{\infty}

\newcommand{\la}{\langle}
\newcommand{\ra}{\rangle}
\newcommand{\wt}{\widetilde}

\newcommand{\n}{\noindent}

\newcommand{\pf}{\noindent{\it Proof.~~}}
\newcommand{\cqd}{\hfill$\Box$}
\newcommand{\be}{\begin{eqnarray*}}
\newcommand{\ee}{\end{eqnarray*}}
\newcommand{\beq}{\begin{equation}}
\newcommand{\eeq}{\end{equation}}

\numberwithin{equation}{section}

\begin{document}



\title{Riesz and Szeg\"o type factorizations for\\
noncommutative Hardy spaces}
\author{Turdebek N. Bekjan and Quanhua Xu}

\date{}

\maketitle


\begin{abstract}

Let $\A$ be a finite subdiagonal algebra in Arveson's sense. Let
$H^p(\A)$ be the associated noncommutative Hardy spaces,
$0<p\le\8$. We extend to the case of all positive indices most
recent results about these spaces, which include notably  the
Riesz, Szeg\"o and inner-outer type factorizations. One new tool
of the paper is the contractivity of the underlying conditional
expectation on $H^p(\A)$ for $p<1$.

\end{abstract}



 \makeatletter
 \renewcommand{\@makefntext}[1]{#1}
 \makeatother \footnotetext{\noindent

\noindent{ T.B.: College of Mathematics and System Sciences,
  Xinjiang
  University, Urumqi 830046 - China\\
  bek@xju.edu.cn\\
\noindent  Q.X.: Laboratoire de
 Math{\'e}matiques,
  Universit{\'e} de Franche-Comt{\'e},
 25030 Besan\c con, cedex - France\\
 qxu@univ-fcomte.fr}
 }


\section{Introduction}


This paper deals with the Riesz and Szeg\"o type factorizations
for noncommutative Hardy spaces associated with a finite
subdiagonal algebra in Arveson's sense \cite{arvesonAJM}.  Let
$\M$ be a finite von Neumann algebra equipped with a normal
faithful tracial state $\t$. Let $\D$ be a von Neumann subalgebra
of $\M$, and let $\F: \M\to\D$ be the unique normal faithful
conditional expectation such that $\t\circ\F=\t$. A {\it finite
subdiagonal algebra} of $\M$ with respect to $\F$ (or $\D$) is a
w*-closed subalgebra $\A$ of $\M$ satisfying the following
conditions
 \begin{enumerate}[i)]
 \item $\A + \A^*$ is w*-dense in $\M$;
 \item $\F$ is multiplicative on $\A$, i.e., $\F(ab)=\F(a)\F(b)$
 for all $a, b\in\A$;
 \item $\A\cap\A^*=\D$.
 \end{enumerate}
We should call the reader's attention to fact that $\A^*$ denotes
in this paper the family of the adjoints of the elements of $\A$,
i.e., $\A^*=\{a^*\;:\; a\in\A\}$. The algebra $\D$ is called the
{\it diagonal} of $\A$. It is proved by Exel \cite{exel} that a
finite subdiagonal algebra $\A$ is automatically {\it maximal} in
the sense that if $\B$ is another subdiagonal algebra with respect
to $\F$ containing $\A$, then $\B=\A$. This maximality yields the
following useful characterization of $\A$:
 \beq\label{maxi-critere}
 \A=\{x\in\M\;:\; \t(xa)=0,\; \forall\; a\in\A_0\},
 \eeq
where $\A_0=\A\cap\ker\F$ (see \cite{arvesonAJM}).

Given $0<p\le\8$ we denote by $L^p(\M)$ the usual noncommutative
$L^p$-space associated with $(\M,\t)$. Recall that $L^\8(\M)=\M$,
equipped with the operator norm.  The norm of $L^p(\M)$ will be
denoted by $\|\cdot\|_p$. For $p<\8$ we define $H^p(\A)$ to be the
closure of $\A$ in $L^p(\M)$, and for $p=\8$ we simply set
$H^\8(\A)=\A$ for convenience. These are the so-called Hardy
spaces associated with $\A$. They are noncommutative extensions of
the classical Hardy spaces on the torus $\mathbb T$. On the other
hand, the theory of matrix-valued analytic functions provides an
important noncommutative example. We refer to \cite{arvesonAJM}
and \cite{px-survey} for more examples.
 We will use the following standard
notation in the theory: If $S$ is a subset of $L^p(\M)$, $[S]_p$
will denote the closure of $S$ in $L^p(\M)$ (with respect to the
w*-topology in the case of $p=\8$). Thus $H^p(\A)=[\A]_p$. Formula
(\ref{maxi-critere}) admits the following $H^p(\A)$ analogue
proved by Saito \cite{saito79}:
 \beq\label{Hp-critere}
 H^p(\A)=\{x\in L^p(\M)\;:\; \t(xa)=0,\; \forall\;
 a\in\A_0\},\quad 1\le p<\8.
 \eeq
Moreover,
 \beq\label{Hp inter Lq}
  H^p(\A)\cap L^q(\M)=H^q(\A), \quad 1\le p<q\le\8.
 \eeq

These noncommutative Hardy spaces have received a lot of attention
since Arveson's pioneer work. We refer the reader notably to the
recent work by Marsalli/West \cite{mars-west} and a series of
newly finished papers by Blecher/Labuschagne \cite{blouter,
blcharac, blbeurling}, whereas more references on previous works
can be found in the survey paper \cite{px-survey}. Most results on
the classical Hardy spaces on the torus have been established in
this noncommutative setting. Here we mention only two of them
directly related with the objective of this paper. The first one
is the Szeg\"o factorization theorem. Already in the fundamental
work \cite{arvesonAJM}, Arveson proved the following factorization
theorem: For any invertible $x\in\M$ there exist a unitary
$u\in\M$ and  $a\in\A$ such that $x=ua$ and $a^{-1}\in\A$. This
theorem is a base of all subsequent works on noncommutative Hardy
spaces. It has been largely improved and extended. The most
general form up to date was newly obtained by Blecher  and
Labuschagne \cite{blouter}: Given $x\in L^p(\M)$ with $1\le
p\le\8$ such that $\De(x)>0$ there exists $h\in H^p(\M)$ such that
$|x|=|h|$. Moreover, $h$ is {\it outer} in the sense that
$[h\A]_p=H^p(\M)$. Here $\De(x)$ denotes the Fuglede-Kadison
determinant of $x$ (see section~2 below for the definition), and
$|x|=(x^*x)^{1/2}$ denotes the absolute value of $x$. We should
emphasize that this result is the (almost) perfect analogue of the
classical Szeg\"o theorem which asserts that given a positive
measurable function $w$ on the torus there exists an outer
function $\f$ such that $w=|\f|$ iff $\log w$ is integrable.

The second result we wish to mention concerns the Riesz
factorization, which asserts that $H^p(\A)=H^q(\A)\cdot H^r(\A)$
for any $1\le p, q,r\le\8$ such that $1/p=1/q+1/r$. More
precisely, given $x\in H^p(\A)$ and $\e>0$ there exist $y\in
H^q(\A)$ and $z\in H^r(\A)$ such that
 $$x=yz\quad\mbox{and}\quad \|y\|_q\,\|z\|_r\le\|x\|_p +\e.$$
This result is proved in \cite{saito79} for $p=q=2$, in
\cite{mars-west} for $r=1$ and independently in \cite{lab-comp}
and in \cite{px-survey} for the general case as above.

Recall that in the case of the classical Hardy spaces the
preceding theorems hold for all positive indices. The problem of
extending these results to the case of indices less than one was
left unsolved in these works. (We mentioned this problem for the
Riesz factorization explicitly in \cite{px-survey}, see the remark
following Theorem~8.3 there). The main purpose of the present
paper is to solve the problem above. As a byproduct, we also
extend all results on outer operators in \cite{blouter} to indices
less than one.

A major obstacle to the solution of the previous problem is the
use of duality, often in a crucial way, in the literature on
noncommutative Hardy spaces.  For instance, duality plays an
important role in proving formulas (\ref{Hp-critere}) and (\ref{Hp
inter Lq}), which are key ingredients for the Riesz factorization
in \cite{mars-west}. In a similar fashion, we will see that their
extensions to indices less than one will be essential for our
proof of the Riesz factorization for all positive indices.

Our key new tool is the contractivity of the conditional
expectation $\F$ on $\A$ with respect to $\|\cdot\|_p$ for
$0<p<1$. Consequently, $\F$ extends to a contractive projection
from $H^p(\A)$ onto $L^p(\D)$. This result is of independent
interest and proved in section~2.

Section~3 is devoted to the Szeg\"o and Riesz type factorizations.
In particular, we extend to all positive indices Marsalli/West's
theorem quoted previously. Section~4 contains some results on
outer operators, notably those in $H^p(\A)$ for $p<1$. This
section can be considered as a complement to the recent work
\cite{blouter}. The last section is devoted to a noncommutative
Szeg\"o formula, which was obtained in \cite{blouter} with the
additional assumption that $\dim\D<\8$.

We will keep all previous notations throughout the paper. In
particular, $\A$ will always denote a finite subdiagonal algebra
of $(\M,\t)$ with diagonal $\D$.


\section{Contractivity of $\F$ on $H^p(\A)$ for $p<1$}


It is well-known that $\F$ extends to a contractive projection
from $L^p(\M)$ onto $L^p(\D)$ for every $1\le p\le\8$.  In
general, $\F$ cannot be, of course, continuously extended to
$L^p(\M)$ for $p<1$. Surprisingly, $\F$ does extend to a
contractive projection on $H^p(\A)$.

\begin{thm}\label{p-contractivity}
 Let $0<p<1$. Then
  \beq\label{p-contractive inequality}
   \forall\; a\in\A\quad \|\F(a)\|_p\le\|a\|_p\,.
  \eeq
 Consequently, $\F$ extends to a contractive projection from
 $H^p(\A)$ onto $L^p(\D)$. The extension will be denoted still by
 $\F$.
  \end{thm}

Inequality (\ref{p-contractive inequality}) is proved by
Labuschagne \cite{lab-szego} for $p=2^{-n}$ and for operators $a$
in $\A$ which are invertible with inverses in $\A$ too.
Labuschagne's proof is a very elegant and simple argument by
induction. It can be adapted to our general situation.

\medskip

\pf Since $\{k2^{-n}\;:\; k, n\in\nat,\; k\ge1\}$ is dense in $(0,
1)$, it suffices to prove  (\ref{p-contractive inequality}) for
$p=k2^{-n}$. Thus we must show
  \beq\label{p-contractive inequality1}
   \forall\; a\in\A\quad \t\big(|\F(a)|^{k2^{-n}}\big)\le
   \t\big(|a|^{k2^{-n}}\big).
  \eeq
This inequality holds for $n=0$ because of the contractivity of
$\F$ on $L^k(\M)$. Now suppose its validity for some $k$ and $n$.
We will prove the same inequality with $n+1$ instead of $n$. To
this end fix $a\in\A$ and $\e>0$. Define, by induction, a sequence
$(x_m)$ by
 $$x_1=(|a|+\e)^{k2^{-n}}\quad\mbox{and}\quad
 x_{m+1}=\frac12\,\big[x_m+(|a|+\e)^{k2^{-n}}\,x_m^{-1}\big].$$
Observe that all $x_m$ belong to the commutative C*-subalgebra
generated by $|a|$. Then it is an easy exercise to show that the
sequence $(x_m)$ is nonincreasing and converges to
$(|a|+\e)^{k2^{-n-1}}$ uniformly (see \cite{lab-szego}). We also
have
  \be
 \t(x_{m+1})
 &=&\frac12\,\big[\t(x_m)+
 \t\big(x_m^{-1/2}(|a|+\e)^{k2^{-n}}\,x_m^{-1/2}\big)\big]\\
 &\ge&\frac12\,\big[\t(x_m)+
 \t\big(x_m^{-1/2}|a|^{k2^{-n}}\,x_m^{-1/2}\big)\big]\\
 &=& \frac12\,\big[\t(x_m)+
 \t\big(|a|^{k2^{-n}}\,x_m^{-1}\big)\big].
 \ee
Now applying Arveson's factorization theorem to each $x_m$, we
find an invertible $b_m\in\A$ with $b_m^{-1}\in\A$ such that
 $$|b_m|=x_m^{2^n/k}\,.$$
Let $p=k2^{-n}$. Then
 \be
 \big\|ab_m^{-1}\big\|_p
 &=&\big\||a|\,b_m^{-1}\big\|_p
 =\big\||a|\,|(b_m^{-1})^*|\big\|_p\\
 &=&\big\||a|\,|b_m|^{-1}\big\|_p
 =\big(\t(|a|^p|b_m|^{-p})\big)^{1/p}\\
 &=&\big(\t(|a|^p\,x_m^{-1})\big)^{1/p}\,,
 \ee
where we have used the commutation between $|a|$ and $|b_m|$ for
the next to the last equality. Therefore, by the induction
hypothesis and the multiplicativity of $\F$ on $\A$
 \be
 \t(x_{m+1})
 &\ge& \frac12\,\big[\t\big(|b_m|^{k2^{-n}}\big)
 +\t\big(|ab_m^{-1}|^{k2^{-n}}\big)\big]\\
 &\ge& \frac12\,\big[\t\big(|\F(b_m)|^{k2^{-n}}\big)+
 \t\big(|\F(a)\F(b_m)^{-1}|^{k2^{-n}}\big)\big].
 \ee
However, by the H\"older inequality
 $$\big(\t\big(|\F(a)|^{k2^{-n-1}}\big)\big)^2
 \le \t\big(|\F(a)\F(b_m)^{-1}|^{k2^{-n}}\big)\,
 \t\big(|\F(b_m)|^{k2^{-n}}\big).$$
It thus follows that
 \be
 \t(x_{m+1})
 &\ge& \frac12\,\big[\t\big(|\F(b_m)|^{k2^{-n}}\big)+
 \big(\t\big(|\F(a)|^{k2^{-n-1}}\big)\big)^2\,
 \big(\t\big(|\F(b_m)|^{k2^{-n}}\big)\big)^{-1}\big]\\
 &\ge& \t\big(|\F(a)|^{k2^{-n-1}}\big).
 \ee
Recalling that $x_m\to (|a|+\e)^{k2^{-n-1}}$ as $m\to\8$, we
deduce
 $$\t\big((|a|+\e)^{k2^{-n-1}}\big)\ge
 \t\big(|\F(a)|^{k2^{-n-1}}\big).$$
Letting $\e\to0$ we obtain inequality (\ref{p-contractive
inequality1}) at the $(n+1)$-th step.\cqd

\begin{cor}\label{multiplicativity of F}
 $\F$ is multiplicative on Hardy spaces. More precisely,
 $\F(ab)=\F(a)\F(b)$ for $a\in H^p(\A)$ and $b\in H^q(\A)$
with $0<p, q\le\8$.
 \end{cor}

\pf  Note that $ab\in H^r(\A)$ for any $a\in H^p(\A)$ and $b\in
H^q(\A)$, where $r$ is determined by $1/r=1/p+1/q$. Thus $\F(ab)$
is well defined. Then the corollary follows immediately from the
multiplicativity of $\F$ on $\A$ and Theorem
\ref{p-contractivity}.\cqd

\medskip

The following is the extension to the case $p<1$ of
Arveson-Labuschagne's Jensen inequality (cf. \cite{arvesonAJM,
lab-szego}). Recall that the Fuglede-Kadison determinant $\De(x)$
of an operator $x\in L^p(\M)$ ($0<p\le\8$) can be defined by
 $$\De(x)=\exp\big(\t(\log|x|)\big)
 =\exp\big(\int_0^\8\log t\,d\nu_{|x|}(t)\big),$$
where $d\nu_{|x|}$ denotes the probability measure on $\real_+$
which is obtained by composing the spectral measure of $|x|$ with
the trace $\t$. It is easy to check that
 $$\De(x)=\lim_{p\to0}\|x\|_p\,.$$
As the usual determinant of matrices, $\De$ is also
multiplicative: $\De(xy)=\De(x)\De(x)$. We refer the reader for
information on determinant to \cite{fug-kad, arvesonAJM} in the
case of bounded operators, and to \cite{brown, haagsch} for
unbounded operators.

\begin{cor}\label{jensen}
 For any $0<p\le\8$ and $x\in H^p(\A)$ we have $\De(\F(x))\le\De(x)$.
 \end{cor}

\pf Let $x\in H^p(\A)$. Then $x\in H^q(\A)$ too for $q\le p$. Thus
by Theorem \ref{p-contractivity}
 $$\|\F(x)\|_q\le\|x\|_q\,.$$
Letting $q\to0$ yields $\De(\F(x))\le\De(x)$.\cqd


\section{Szeg\"o and Riesz factorizations}


The following result is a Szeg\"o type factorization theorem. It
is stated in \cite{px-survey} without proof (see the remark
following Theorem~8.1 there). We take this opportunity to provide
a proof. It is an improvement of the previous factorization
theorems of Arveson \cite{arvesonAJM} and Saito \cite{saito79}. As
already quoted in the introduction, Blecher and Labuschagne newly
obtained a Szeg\"o factorization for any $w\in L^p(\M)$ with $1\le
p\le\8$ such that $\De(w)>0$ (see the next section for more
details). Note that the property that $h^{-1}\in H^q(\A)$ whenever
$w^{-1}\in L^q(\M)$ will be important for our proof of the Riesz
factorization below. Let us also point out that although not in
full generality, this result has hitherto been strong enough for
applications in the literature. See Theorem \ref{szego bis} below
for an improvement.

\begin{thm}\label{szego}
 Let $0<p, q\le\8$. Let $w\in L^p(\M)$ be an invertible operator
such that $w^{-1}\in L^q(\M)$. Then there exist a unitary $u\in\M$
and $h\in H^p(\A)$  such that $w=uh$ and $h^{-1}\in H^q(\A)$.
 \end{thm}

\pf  We first consider the case $p=q=2$. The proof of this special
case is modelled on Arveson's original proof of his Szeg\"o
factorization theorem (see also \cite{saito79}). Let $x$ be the
orthogonal projection of $w$ in $[w\A_0]_2$, and set $y=w-x$. Thus
$y\perp[w\A_0]_2$; whence $y\perp[y\A_0]_2$. It follows that
 $$\forall\; a\in\A_0\quad \t(y^*ya)=0.$$
Hence by (\ref{Hp-critere}), $y^*y\in H^1(\A)=[\A]_1$, and
$y^*y\in[\A^*]_1$ too. On the other hand, it is easy to see that
$[\A]_1\cap [\A^*]_1=L^1(\D)$. Indeed, if $a\in [\A]_1\cap
[\A^*]_1$, then $\t(ab)=0$ for any $b\in \A_0+\A_0^*$; so
$\t(ab)=\t(\F(a)b)$ for any $b\in \A+\A^*$. It follows that
$a=\F(a)\in L^1(\D)$. Consequently, $y^*y\in L^1(\D)$, so $|y|\in
L^2(\D)$.

Regarding $\M$ as a von Neumann algebra acting on $L^2(\M)$ by
left multiplication, we claim that $y$ is cyclic for $\M$. This is
equivalent to showing that $y$ is separating for the commutant of
$\M$. However, this commutant coincides with the algebra of all
right multiplications on $L^2(\M)$ by the elements of $\M$. Thus
we are reduced to prove that if $z\in\M$ is such that $yz=0$, then
$z=0$. We have:
 $$0=\t(z^*y^*yz)=\t(|y|^2\,|z^*|^2)=\t(|y|^2\,\F(|z^*|^2))
 =\|yd\|_2^2\,,$$
where $d=\F(|z^*|^2)^{1/2}\in\D$; whence $yd=0$. Choose a sequence
$(a_n)\subset\A_0$ such that
 \beq\label{x limit} x=\lim wa_n.\eeq
 Then (recalling that $w^{-1}\in L^2(\M)$)
 $$0=\t(w^{-1}yd)=\lim_n\t(w^{-1}(w- wa_n)d)
 =\t(d)-\lim_n \t(a_nd)=\t(d).$$
It follows that $d=0$, so by virtue of the faithfulness of $\F$,
$z=0$ too. This yields our claim. Therefore, $[\M y]_2=L^2(\M)$.
It turns out that the right support of $y$ is $1$. Since $\M$ is
finite, the left support of $y$ is also equal to $1$, so $y$ is of
full support. Consequently, $[y\M]_2=L^2(\M)$ too.

Let $y=u|y|$ be the polar decomposition of $y$. Then $u$ is a
unitary in $\M$. Let $h=u^*w$. We are going to prove that $h\in
H^2(\A)$. To this end we first note the following orthogonal
decomposition of $L^2(\M)$:
 \beq\label{dec}
 L^2(\M)= [y\A_0]_2\oplus [y\D]_2\oplus[y\A_0^*]_2\,.
 \eeq
Indeed, for any $a\in \A$ and $b\in\A_0$ we have
 $$\la ya,\; yb^*\ra=\t(by^*ya)=\t(|y|^2ab)=0;$$
so $[y\A_0]_2\oplus [y\D]_2\oplus[y\A_0^*]_2$ is really an
orthogonal sum. On the other hand, by the previous paragraph, we
see that
 $$L^2(\M)=[y\M]_2 \subset
 [y\A_0]_2\oplus [y\D]_2\oplus[y\A_0^*]_2\,.$$
 Therefore,  decomposition (\ref{dec}) follows.  Applying $u^*$
to both sides of (\ref{dec}), we deduce
  \be
  L^2(\M)
  &=& [u^*y\A_0]_2\oplus [u^*y\D]_2\oplus[u^*y\A_0^*]_2\\
  &=&[|y|\A_0]_2\oplus [|y|\D]_2\oplus[|y|\A_0^*]_2\,.
  \ee
 Since $|y|\in L^2(\D)$, $[|y|\A_0]_2\subset [A_0]_2$, and
 similarly for the two other terms on the right. Therefore,
 \be
  L^2(\M)
  &=&[|y|\A_0]_2\oplus [|y|\D]_2\oplus[|y|\A_0^*]_2\\
  &\subset& [\A_0]_2\oplus [\D]_2\oplus[\A_0^*]_2=L^2(\M)\,.
  \ee
 Hence
  \beq\label{dec1}
  [|y|\A_0]_2=[\A_0]_2,\quad [|y|\D]_2=[\D]_2, \quad
  [|y|\A_0^*]_2=[\A_0^*]_2\,.
  \eeq
Passing to adjoints, we also have
 $$[\A_0|y|]_2=[\A_0]_2,\quad [\D|y|]_2=[\D]_2, \quad
  [\A_0^*|y|]_2=[\A_0^*]_2\,.$$
Now it is easy to show that $h=u^*w\in H^2(\A)$. Indeed, since
$y\perp[w\A_0]$, $\t(y^* wa)=0$ for all $a\in\A_0$; so $\t(a|y|
u^*w)=0$. However, $[\A_0|y|]_2=[\A_0]_2$. Thus
 $$\forall\; a\in H^2_0(\A)\quad \t(a h)=0.$$
Hence by (\ref{Hp-critere}),  $h\in H^2(\A)$, as desired.

It remains to show that $h^{-1}\in H^2(\A)$. To this end we first
observe that $\F(h)\F(h^{-1})=1$. Indeed, given $d\in\D$ we have,
by (\ref{x limit})
 \be
 \t\big(\F(h)\F(h^{-1})|y|d\big)
 &=&\t\big(h^{-1}|y|d\F(h)\big)
 =\t\big(w^{-1}u|y|d\F(h)\big)\\
 &=&\lim_n\t\big(w^{-1}(w-wa_n)d\F(h)\big)
 =\t\big(d\F(h)\big)\\
 &=&\t(hd)
 =\t(u^*wd)=\t(u^*yd)=\t(|y|d),
 \ee
where we have used the fact that
 $$\t(u^*xd)=\lim_n\t(u^*wa_nd)=\lim_n\t(ha_nd)=0.$$
Since $[|y|\D]_2=L^2(\D)$, we deduce our observation. Therefore,
$\F(h)$ is invertible and its inverse is $\F(h^{-1})$. On the
other hand, by (\ref{x limit})
 $$\F(h)=\lim_n\F(u^*(y+wa_n))
 =\F(|y|)+ \lim_n\F(ha_n)=u^*y.$$
Hence,
 $$u=y\F(h)^{-1}=y\F(h^{-1}).$$
Now let $a\in \A_0$. Then
 \be
 \t(h^{-1}a)
 =\t(w^{-1}ua)=\t\big(w^{-1}y\F(h^{-1})a\big)
 =\lim_n\t\big(w^{-1}(w-wa_n)\F(h^{-1})a\big)
 =0.
 \ee
It follows that $h^{-1}\in H^2(\A)$. Therefore, we are done in the
case $p=q=2$.

The general case can be easily reduced to this special one.
Indeed, if $p\ge2$ and $q\ge2$, then given $w\in L^p(\M)$ with
$w^{-1}\in L^q(\M)$, we can apply the preceding part and then find
a unitary $u\in\M$ and $h\in H^2(\A)$ such that $w=uh$ and
$h^{-1}\in H^2(\A)$. Then $h=u^*w\in L^p(\M)$, so $w\in
H^2(\A)\cap L^p(\M)=H^p(\A)$ by (\ref{Hp inter Lq}). Similarly,
$h^{-1}\in H^q(\A)$.

Suppose $\min(p,\; q)<2$. Choose an integer $n$ such that
$\min(np,\; nq)\ge2$. Let $w=v|w|$ be the polar decomposition of
$w$. Note that $v\in\M$ is a unitary. Write
 $$w=v|w|^{1/n}\, |w|^{1/n}\,\cdots\, |w|^{1/n}
 =w_1w_2\cdots w_n,$$
where $w_1=v|w|^{1/n}$ and $w_k=|w|^{1/n}$ for $2\le k\le n$.
Since $w_k\in L^{np}(\M)$ and $w_k^{-1}\in L^{nq}(\M)$, by what is
already proved we have a factorization
 $$w_n=u_n h_n$$
with $u_n\in\M$ a unitary, $h_n\in H^{np}(\A)$ such that
$h_n^{-1}\in H^{nq}(\A)$. Repeating this argument, we again get a
same factorization for $w_{n-1}u_n$:
 $$w_{n-1}u_n=u_{n-1}h_{n-1}\,;$$
and then for $w_{n-2}u_{n-1}$, and so on. In this way, we obtain a
factorization:
 $$w=u h_1\,\cdots\, h_n,$$
where $u\in\M$ is a unitary, $h_k\in H^{np}(\A)$ such that
$h_k^{-1}\in H^{nq}(\A)$. Setting $h=h_1\,\cdots\, h_n$, we then
see that $w=uh$ is the desired factorization. Hence the proof of
the theorem is complete.\cqd

\begin{rk}\label{szegork}
 Let $w\in L^2(\M)$ be an invertible operator such that
$w^{-1}\in L^2(\M)$. Let $w=uh$ be the factorization in Theorem
\ref{szego}. The preceding proof shows that $[h\A]_2=H^2(\A)$.
Indeed, it is clear that $[y\A]_2\subset[w\A]_2$. Using
decomposition (\ref{dec}), we get
 $$[w\A]_2\ominus [y\A]_2=[w\A]_2\cap [y\A_0^*]_2\,.$$
Now for any $a\in\A$ and $b\in\A_0$,
 $$\la wa,\; yb^*\ra=\t(y^*wab)=0$$
since $y\perp [w\A_0]$. It then follows that $[w\A]_2\ominus
[y\A]_2=\{0\}$, so $[w\A]_2=[y\A]_2$. Hence, by (\ref{dec1})
 $$[h\A]_2=[u^*w\A]_2 =[u^*y\A]_2=[|y|\A]_2=H^2(\A).$$
 \end{rk}

We turn to the Riesz factorization. We  first need to extend
(\ref{Hp inter Lq}) to all indices.

\begin{prop}\label{intersection}
 Let $0<p<q\le\8$. Then
  $$H^p(\A)\cap L^q(\M)=H^q(\A)\quad\mbox{and}\quad
  H^p_0(\A)\cap L^q(\M)=H^q_0(\A),$$
where $H^p_0(\A)=[\A_0]_p$.
 \end{prop}

\pf  It is obvious that $H^q(\A)\subset H^p(\A)\cap L^q(\M)$. To
prove the converse inclusion, we first consider the case $q=\8$.
Thus let $x\in H^p(\A)\cap \M$. Then by Corollary
\ref{multiplicativity of F},
 $$\forall\; a\in\A_0\quad \F(xa)=\F(x)\F(a)=0.$$
Hence by (\ref{maxi-critere}), $x\in \A$.

Now consider the general case. Fix an $x\in H^p(\A)\cap L^q(\M)$.
Applying Theorem \ref{szego} to $w=(x^*x+1)^{1/2}$, we get an
invertible $h\in H^q(\A)$ such that
 $$h^*h=x^*x+1\quad\mbox{and}\quad h^{-1}\in\A.$$
Since $h^*h\le x^*x$, there exists a contraction $v\in\M$ such
that $x=vh$. Then $v=xh^{-1}\in H^p(\A)\cap\M$, so  $v\in\A$.
Consequently, $x\in \A\cdot H^q(\A)=H^q(\A)$. Thus we proved the
first equality. The second is then an easy consequence. For this
it suffices to note that $H^p_0(\A)=\{x\in H^p(\A)\;:\;\F(x)=0\}$.
The later equality follows from the continuity of $\F$ on
$H^p(\A)$. \cqd

\begin{thm}\label{riesz}
 Let $0< p, q, r\le\8$ such that $1/p=1/q+1/r$.
Then for $x\in H^p(\A)$ and $\e>0$ there exist $y\in H^q(\A)$ and
$z\in H^r(\A)$ such that
 $$x=yz\quad\mbox{and}\quad \|y\|_q\,\|z\|_r\le \|x\|_p+\e.$$
Consequently,
 $$\|x\|_p=\inf\big\{\|y\|_q\,\|z\|_r\;:\; x=yz,\; y\in H^q(\A),\;
 z\in H^r(\A) \big\}.$$
\end{thm}

\pf The case where $\max(q,\;r)=\8$ is trivial. Thus we assume
both $q$ and $r$ to be finite. Let $w=(x^*x+\e)^{1/2}$. Then $w\in
L^p(\M)$ and $w^{-1}\in\M$. Let $v\in\M$ be a contraction  such
that $x=vw$. Now applying Theorem \ref{szego} to $w^{p/r}$, we
have: $w^{p/r}=uz$, where $u$ is a unitary in $\M$ and $z\in
H^r(\A)$ such that $z^{-1}\in\A$. Set $y=vw^{p/q}\,u$. Then
$x=yz$, so $y=xz^{-1}$. Since $x\in H^p(\A)$ and $z^{-1}\in\A$,
$y\in H^p(\A)$. On the other hand, $y$ belongs to $L^{q}(\M)$ too.
Therefore, $y\in H^q(\A)$ by virtue of Proposition
\ref{intersection}. The norm estimate is clear.\cqd

\begin{rk}
 It is unknown at the time of this writing whether the
infimum in Theorem \ref{riesz} is attained. We will see in
section~4 that the answer is affirmative if additionally
$\De(x)>0$.
\end{rk}


\section{Outer operators}


We consider in this section outer operators.  All results below on
the left and right outers are due to Blecher and Labuschagne
\cite{blouter} in the case of indices not less than one. The
notion of bilaterally outer is new. We start with the following
result.

\begin{prop}\label{outer: independence}
 Let $0<p<q\le\8$ and let $h\in H^q(\A)$.
Then
 \begin{enumerate}[\rm i)]
 \item $[h\A]_p=H^p(\A)$ iff $[h\A]_q=H^q(\A)$;
 \item $[\A h]_p=H^p(\A)$ iff $[\A h]_q=H^q(\A)$;
 \item $[\A h\A]_p=H^p(\A)$ iff $[\A h\A]_q=H^q(\A)$.
 \end{enumerate}
 \end{prop}

\pf We prove only the third equivalence. The proofs of the two
others are similar (and even simpler). It is clear that $[\A
h\A]_q=H^q(\A)\;\Rightarrow\; [\A h\A]_p=H^p(\A)$. To prove the
converse implication we first consider the case $q\ge1$. Let $q'$
be the conjugate index of $q$. Let $x\in L^{q'}(\M)$ be such that
 $$\forall\; a, b\in\A\quad \t(xahb)=0.$$
Then $xah\in H_0^1(\A)$ for any $a\in\A$ by virtue of
(\ref{Hp-critere}) (more rigorously, its $H^p_0$-analogue as in
Proposition \ref{intersection}). On the other hand, by the
assumption that $[\A h\A]_p=H^p(\A)$, there exist two sequences
$(a_n), (b_n)\subset \A$ such that
 $$\lim_na_nhb_n=1\quad \mbox{in}\quad H^p(\A).$$
Consequently,
 $$\lim_nxa_nhb_n=x\quad \mbox{in}\quad L^r(\M),$$
where $1/r=1/q' + 1/p$. Since $xa_nhb_n=(xa_nh)b_n\in
H^1_0(\A)\subset H^r_0(\A)$, we deduce that $x\in H^r_0(\A)$.
Therefore, $x\in H^r_0(\A)\cap L^{q'}(\M)$, so by Proposition
\ref{intersection}, $x\in H_0^{q'}(\A)$. Hence, $\t(xy)=0$ for all
$y\in H^q(\A)$. Thus $[\A hA]_q=H^q(\A)$.

Now assume $q<1$. Choose an integer $n$ such that $np\ge2$. By the
proof of Theorem \ref{riesz} and Remark \ref{szegork}, we deduce a
factorization:
 $$h=h_1\,h_2\,\cdots\,h_n\,,$$
where $h_k\in H^{np}(\A)$ for every $1\le k\le n$ and
$[h_k\A]_2=H^2(\A)$ for $2\le k\le n$.  By the left version (i.e.,
part i)) of the previous case already proved,  we also have
$[h_k\A]_{np}=H^{np}(\A)$ and $[h_k\A]_{nq}=H^{nq}(\A)$ for $2\le
k\le n$. Let us deal with the first factor $h_1$. Using $[\A
h\A]_p=H^p(\A)$ and $[h_k\A]_{np}=H^{np}(\A)$  for $2\le k\le n$,
we see that $[\A h_1\A]_p=H^p(\A)$; so again $[\A
h_1\A]_{nq}=H^{nq}(\A)$ by virtue of the first part. It is then
clear that $[\A h\A]_q=H^q(\A)$. \cqd

\medskip

The previous result justifies the relative independence of the
index $p$ in the following definition.

\begin{definition}
 Let $0<p\le\8$. An operator $h\in H^p(\A)$ is called {\it left outer},
{\it right outer} or {\it bilaterally outer} according to
$[h\A]_p=H^p(\A)$, $[\A h]_p=H^p(\A)$ or $[\A h\A]_p=H^p(\A)$.
 \end{definition}

\begin{rk}\label{rk outer}
 It is easy to see that if $h$ is left outer or right outer, $h$
is of full support (i.e., $h$ is injective and of dense range).
There exist, however, bilaterally outers which are not of full
support. For example, consider the case where $\A=\M=\mathbb M_n$,
the full algebra of $n\times n$ complex matrices, equipped with
the normalized trace. Then every $e_{ij}$ is bilaterally outer,
where the $e_{ij}$ are the canonical matrix units of $\mathbb
M_n$. A less trivial case is the following. Let $\T$ be the unit
circle equipped with normalized Haar measure. Let $\M=L^\8(\mathbb
T)\bar\otimes\mathbb M_n=L^\8(\mathbb T;\mathbb M_n)$, and let
$\A=H^\8(\mathbb T;\mathbb M_n)$, the algebra of $\mathbb
M_n$-valued bounded analytic functions in  the unit disc of the
complex plane. Let $\f\in H^p(\T)$ be an outer function. Then
$h=\f\otimes e_{ij}$ is bilaterally outer with respect to $\A$.
 \end{rk}

\begin{thm}\label{outer-critere}
 Let $0<p\le\8$ and $h\in H^p(\A)$.
 \begin{enumerate}[\rm i)]
 \item If $h$ is left or right outer, then
$\De(h)=\De(\F(h))$. Conversely, if $\De(h)=\De(\F(h))$ and
$\De(h)>0$, then $h$ is left and right outer $($so bilaterally
outer too$)$.
 \item If $\A$ is antisymmetric $($i.e., $\dim\D=1)$
and $h$ is bilaterally outer, then $\De(h)=\De(\F(h))$.
 \end{enumerate}
 \end{thm}

\pf i) This part is proved in \cite{blouter} for $p\ge1$. Assume
$h$ is left outer. Let $d\in\D$. Using Theorem
\ref{p-contractivity}, we obtain
 $$\|\F(h)d\|_p=\inf\big\{\|hd+x_0\|_p\;:\; x\in H^p_0(\A)\big\}.$$
On the other hand,
 $$[h\A_0]_p=\big[[h\A]_p\A_0\big]_p
 =\big[[\A]_p\A_0\big]_p=[\A_0]_p=H^p_0(\A).$$
Therefore,
 $$\|\F(h)d\|_p=\inf\big\{\|h(d+a_0)\|_p\;:\; a_0\in \A_0\big\}.$$
Recall the following characterization of $\De(x)$ from
\cite{blouter}:
 \beq\label{det formula}
 \De(x)=\inf\big\{\|xa\|_p\;:\; a\in\A,\;
 \De(\F(a))\ge1\big\}.
 \eeq
Now using this formula twice, we obtain
 \be
 \De(\F(h))
 &=&\inf\big\{\|\F(h)d\|_p\;:\; d\in\D, \De(d)\ge1\big\}\\
 &=&\inf\big\{\|h(d+a_0)\|_p\;:\; d\in\D, \De(d)\ge1,\;a_0\in \A_0\big\}\\
 &=&\De(h).
 \ee

Let us show the converse under the additional assumption that
$\De(h)>0$. We will use the case $p\ge1$ already proved in
\cite{blouter}. Thus assume $p<1$. Choose an integer $n$ such that
$np\ge1$.  By Theorem \ref{riesz}, there exist $h_1,\,...\,,
h_n\in H^{np}(\A)$ such that $h=h_1\,\cdots\, h_n$. Then
 $\De(h) =\De(h_1)\,\cdots\,\De(h_n);$
so $\De(h_k)>0$ for all $1\le k\le n$. On the other hand,  by
Arveson-Labuschagne's Jensen inequality \cite{arvesonAJM,
lab-szego} (or Corollary \ref{jensen}), $\De(\F(h_k))\le
\De(h_k)$. However,
 \be
 \De(\F(h))=\De(\F(h_1))\,\cdots\,\De(\F(h_n))
 \le \De(h_1)\,\cdots\,\De(h_n)
 =\De(h)=\De(\F(h)).
 \ee
It then follows that $\De(\F(h_k))= \De(h_k)$ for all $k$. Now
$h_k\in H^{np}(\A)$ with $np\ge1$, so  $h_k$ is left and right
outer. Consequently, $h$ is left and right outer.

ii) This proof is similar to that of the first part of  i). We
will use the following variant of (\ref{det formula})
 \beq\label{det formula bis}
 \De(x)=\inf\big\{\|axb\|_p\;:\; a, b\in\A,\;
 \De(\F(a))\ge1,\;\De(\F(b))\ge1\big\}
 \eeq
for every $x\in L^p(\M)$. This formula immediately follows from
(\ref{det formula}). Indeed, by (\ref{det formula}) and the
multiplicativity of $\De$
 \be
 &&\inf\big\{\|axb\|_p\;:\; a,
 b\in\A,\;\De(\F(a))\ge1,\;\De(\F(b))\ge1\big\}\\
 &&=\inf\big\{\De(ax)\;:\; a\in\A,\;\De(\F(a))\ge1\big\}\\
 &&=\inf\big\{\De(a)\De(x)\;:\; a\in\A,\;\De(\F(a))\ge1\big\}
 =\De(x).
 \ee
Now assume $h\in H^p(\A)$ is bilaterally outer and $\A$ is
antisymmetric. Then $\F(h)$ is a multiple of the unit of $\M$. As
in the proof of i), We have
  \begin{eqnarray}\label{det formula2}
  \begin{array}{ccl}
 \begin{displaystyle}\|\F(h)\|_p\end{displaystyle}
 &=&\begin{displaystyle}
 \inf\big\{\|h+x\|_p\;:\; x\in H^p_0(\A)\big\}
 \end{displaystyle}\\
 &=&\begin{displaystyle}
 \inf\big\{\|h+ahb_0\|_p\;:\; a\in\A,\; b_0\in\A_0\big\}
 \end{displaystyle}.
 \end{array}
 \end{eqnarray}
Using $\dim\D=1$, we easily check that
 \beq\label{det formula3}
 \inf\big\{\|h+ahb_0\|_p\;:\; a\in\A,\; b_0\in\A_0\big\}
 =\inf\big\{\|(1+a_0)h(1+b_0)\|_p\;:\; a_0, b_0\in\A_0\big\}.
 \eeq
Indeed, it suffices to show that both sets $\{h+ahb_0:\, a\in\A,\;
b_0\in\A_0\}$ and $\{(1+a_0)h(1+b_0):\, a_0, b_0\in\A_0\}$ are
dense in $\{x\in H^p(\A):\,\F(x)=\F(h)\}$. The first density
immediately follows from the density of $\A h\A_0$ in $H^p_0(\A)$.
On the other hand, let $x\in H^p(\A)$ with $\F(x)=\F(h)$ and let
$a_n, b_n\in\A$ such that
 $\lim_n a_nhb_n=x.$
By Theorem \ref{p-contractivity},
 $$\lim_n \F(a_n)\F(h)\F(b_n)=\F(x).$$
Since $\F(x)=\t(x)1=\t(h)1=\F(h)\neq0$, we deduce that $\lim_n
\t(a_n)\t(b_n)=1.$ Thus replacing $a_n$ and $b_n$ by $a_n/\t(a_n)$
and $b_n/\t(b_n)$, respectively, we can assume that $a_n=1+\wt
a_n$ and $b_n=1+\wt b_n$ with $\wt a_n, \wt b_n\in\A_0$; whence
the desired density of $\{(1+a_0)h(1+b_0):\, a_0, b_0\in\A_0\}$ in
$\{x\in H^p(\A):\,\F(x)=\F(h)\}$. Finally, combining  (\ref{det
formula bis}), (\ref{det formula2}) and (\ref{det formula3}), we
get $\De(\F(h))=\De(h)$. \cqd

\begin{rk}\label{rk bi-outer}
 The assumption that $\A$ is antisymmetric in Theorem
\ref{outer-critere}, ii) cannot be removed in general, as shown by
the following example. Keep the notation introduced in Remark
\ref{rk outer} and consider the case where $\M=L^\8(\T; \mathbb
M_2)$ and $\A=H^\8(\T; \mathbb M_2)$. Let $\f_1$ and $\f_2$ be two
outer functions in $H^p(\T)$, and let $h=\f_1\otimes e_{11} +
z\f_2\otimes e_{22}$, where $z$ denotes the identity function on
$\T$. Then it is easy to check that $h$ is bilaterally outer and
 $$\De(h)=\exp\Big(\frac12 \int_{\T}\log|\f_1| +
 \frac12 \int_{\T}\log|\f_2|\Big)>0.$$
However, $\F(h)=\f_1(0) e_{11}$, so $\De(\F(h))=0$.
 \end{rk}

The following is an immediate consequence of Theorem
\ref{outer-critere}. We do not know, however, whether the
condition $\De(h)>0$ in i) can be removed or not.

\begin{cor}
 Let $h\in H^p(\A)$, $0<p\le\8$.
 \begin{enumerate}[\rm i)]
 \item  If $\De(h)>0$, then $h$ is left outer iff
$h$ is right outer.
 \item Assume that $\A$ is antisymmetric.
Then the following properties are equivalent:
 \begin{enumerate}[$\bullet$]
 \item $h$ is left outer;
 \item $h$ is right outer;
 \item $h$ is bilaterally outer;
 \item $\De(\F(h))=\De(h)>0$.
 \end{enumerate}
 \end{enumerate}
 \end{cor}

We  will say that $h$ is {\it outer} if it is at the same time
left and right outer. Thus if $h\in H^p(\A)$ with $\De(h)>0$, then
$h$ is outer iff $\De(h)=\De(\F(h))$. Also in the case where $\A$
is antisymmetric, an $h$ with $\De(h)>0$ is outer iff it is left,
right or bilaterally outer.

\begin{cor}
 Let $h\in H^p(\A)$ such that $h^{-1}\in H^q(\A)$ with
 $0<p,q\le\8$. Then $h$ is outer.
 \end{cor}

\pf By the multiplicativity of $\De$, $\De(h)\De(h^{-1})=1$ and
$\De(\F(h))\De(\F(h^{-1}))=1$. Thus by Jensen's inequality
(Corollary \ref{jensen}),
 $$\De(h)=\De(h^{-1})^{-1}\le \De(\F(h^{-1}))^{-1}=\De(\F(h));$$
whence the assertion because of Theorem  \ref{outer-critere}.\cqd

\medskip

The following improves Theorem \ref{szego}.

\begin{thm}\label{szego bis}
 Let $w\in L^p(\M)$ with $0<p\le\8$ such that $\De(w)>0$. Then
there exist a unitary $u\in\M$ and an outer $h\in H^p(\A)$ such
that $w=uh$.
 \end{thm}

\pf Based on the case  $p\ge1$ from \cite{blouter}, the proof
below is similar to the end of the proof of Theorem \ref{szego}.
For simplicity we consider only the case where $p\ge 1/2$. Write
the polar decomposition of $w$:  $w=v|w|$. Applying \cite{blouter}
to $|w|^{1/2}$ we get a factorization: $|w|^{1/2}=u_2h_2$ with
$u_2$ unitary and $h_2\in H^{2p}(\A)$ left outer. Since
$\De(h_2)>0$, $h_2$ is also right outer; so $h_2$ is outer.
Similarly, we have: $v|w|^{1/2}u_2=u_1h_1$. Then $u=u_1$ and
$h=h_1h_2$ yield the desired factorization of $w$.\cqd

\medskip

The following is the inner-outer factorization for operators in
$H^p(\A)$, which is already in \cite{blouter} for $p\ge1$.

\begin{cor}\label{i-o}
 Let $0<p\le\8$ and $x\in H^p(\A)$ with $\De(x)>0$. Then there
 exist a unitary $u\in\A$ $($inner$)$ and an outer $h\in H^p(\A)$
 such that $x=uh$.
 \end{cor}

\pf Applying the previous theorem, we get $x=uh$ with $h$ outer
and $u$ a unitary in $\M$. Let $a_n\in\A$ such that $\lim ha_n=1$
in $H^p(\A)$. Then $u=\lim xa_n$ in $H^p(\A)$ too; so $u\in
H^p(\A)\cap\M$. By Proposition \ref{intersection}, $u\in\A$.\cqd

\begin{rk}
 The condition $\De(x)>0$ cannot be removed in general.
Indeed, if $h$ is outer, then $h$ is of full support (see Remark
\ref{rk bi-outer}). It follows that $x$ is of full support too if
$x$ admits an inner-outer factorization as above. Consider, for
instance, the example in Remark \ref{rk bi-outer}. Then for any
$\f\in H^p(\mathbb T)$ the operator $x=\f\ot e_{11}\in H^p(\A)$ is
not of full support.
 \end{rk}

\begin{cor}
 Let $0<p\le\8$ and $h\in H^p(\A)$ with $\De(h)>0$. Then $h$ is
outer iff for any $x\in H^p(\A)$ with $|x|=|h|$ we have
$\De(\F(x))\le \De(\F(h))$.
 \end{cor}

\pf Assume $h$ outer. Then by Corollary \ref{jensen} and Theorem
\ref{outer-critere},
 $$\De(\F(x))\le \De(x)=\De(h)=\De(\F(h)).$$
Conversely, let $h=uk$ be the decomposition given by Theorem
\ref{szego bis} with $k$ outer. Then
 $$\De(h)=\De(k)=\De(\F(k))\le \De(\F(h));$$
so $h$ is outer by Theorem \ref{outer-critere}.\cqd

\begin{cor}
 Let $0< p, q, r\le\8$ such that $1/p=1/q+1/r$. Let $x\in H^p(\A)$
be such that $\De(x)>0$. Then there exist $y\in H^q(\A)$ and $z\in
H^r(\A)$ such that
 $$x=yz\quad\mbox{and}\quad \|x\|_p=\|y\|_q\,\|z\|_r\,.$$
\end{cor}

\pf This proof is similar to that of Theorem \ref{riesz}. Instead
of Theorem \ref{szego}, we now use Theorem \ref{szego bis}.
Indeed, by the later theorem, we can find a unitary $u_2\in\M$ and
an outer $h_2\in H^{p/r}(\A)$ such that $|x|^{p/r}=u_2h_2$. Once
more applying this theorem to $v|x|^{p/q}\,u_2$, we have a similar
factorization: $v|x|^{p/q}\,u_2=u_1h_1$, where $v$ is the unitary
in the polar decomposition of $x$. Since $h_1$ and $h_2$ are
outer, we deduce, as in the proof of Corollary \ref{i-o}, that
$u_1\in\A$. Then $y=u_1h_1$ and $z=h_2$ give the desired
factorization of $x$.\cqd


\section{A noncommutative Szeg\"o formula}


Let $w\in L^1(\mathbb T)$ be a positive function and let
$d\mu=wdm$. Then we have the following well-known Szeg\"o formula
\cite{szego}:
 $$\inf\big\{\int_{\mathbb T}|1-f|^2d\mu\;:\;
 f \;\mbox{mean zero analytic polynomial}\big\}=\exp\int_{\mathbb
 T}\log wdm.$$
This formula was later proved for any positive measure $\mu$ on
$\mathbb T$ independently by  Kolmogorov/Krein \cite{kol} and
Verblunsky \cite{ver}. Then the singular part of $\mu$ with
respect to the Lebesgue measure $dm$ does not contribute to the
preceding infimum  and $w$ on the right hand side is the density
of the absolute part of $\mu$ (also see \cite{garn}). This latter
result was extended to the noncommutative setting in
\cite{blouter}. More precisely, let $\o$ be a positive linear
functional on $\M$, and let $\o=\o_n+\o_s$ be the decomposition of
$\o$ into its normal and singular parts. Let $w\in L^1(\M)$ be the
density of $\o_n$ with respect to $\t$, i.e., $\o_n=\t(w\,\cdot)$.
Then Blecher and Labuschagne proved that if $\dim\D<\8$,
 $$\De(w)=\inf\big\{\o(|a|^2)\;:\; a\in\A,\; \De(\F(a))\ge1\big\}.$$
It is left open in \cite{blouter} whether the condition
$\dim\D<\8$ can be removed or not. We will solve this problem in
the affirmative. At the same time, we show that the square in the
above formula can be replaced by any power $p$.

\begin{thm}\label{szego formula}
 Let $\o=\o_n+\o_s$ be as above and $0<p<\8$. Then
 $$\De(w)=\inf\big\{\o(|a|^p)\;:\; a\in\A,\; \De(\F(a))\ge1\big\}.$$
 \end{thm}

\pf Let
 $$\d(\o)=\inf\{\o(|a|^p)\;:\; a\in\A, \;\De(\F(a))\ge1\}.$$
First we show that
 $$\d(\o)=\inf\{\o(x)\;:\; x\in\M_+^{-1},\; \De(x)\ge1\},$$
where $\M_+^{-1}$ denotes the family of invertible positive
operators in $\M$ with bounded inverses. Given any
$x\in\M_+^{-1}$, by Arveson's factorization theorem there exists
$a\in\A$ such that $|a|=x^{1/p}$ and $a^{-1}\in\A$. Then
$x=|a|^p$, so $\De(x)=\De(|a|^p)=\De(a)^p$. Since $a$ is
invertible with $a^{-1}\in\A$, by Jensen's formula in
\cite{arvesonAJM}, $\De(a)=\De(\F(a))$. It then follows that
 $$\d(\o)\le\inf\{\o(x)\;:\; x\in\M_+^{-1},\; \De(x)\ge1\}.$$
The converse inequality is easier. Indeed, given $a\in\A$ with
$\De(\F(a))\ge1$ and $\e>0$, set $x=|a|^p+\e$. Then
$x\in\M_+^{-1}$ and $\De(x)\ge\De(a)^p\ge \De(\F(a))^p$ by virtue
of Jensen's inequality. Since
$\lim_{\e\to0}\o(|a|^p+\e)=\o(|a|^p)$, we deduce the desired
converse inequality.

Next we show that $\d(\o)=\d(\o_n)$.  The singularity of $\o_s$
implies that there exists  an increasing net $(e_i)$ of
projections in $\M$ such that $e_i\to 1$ strongly and
$\o_s(e_i)=0$ for every $i$ (see \cite[III.3.8]{tak-I}). Let
$\e>0$. Set
 $$x_i=\e^{\t(e_i)-1}(e_i +\e e_i^\perp), \quad \mbox{where}
 \quad e^\perp=1-e.$$
Clearly, $x_i\in \M_+^{-1}$ and $\De(x_i)=1$. Let $x\in\M_+^{-1}$
and $\De(x)\ge1$. Then $\De(x_ixx_i)=\De(x)\ge1$, and $x_ixx_i\to
x$ in the w*-topology. On the other hand, note that
 $$\o_s(x_ixx_i)=\e^{2\t(e_i)}\o_s(e_i^\perp xe_i^\perp).$$
Therefore,
 \be
 \d(\o)
 &\le& \limsup\o(x_ixx_i)=\o_n(x)+ \limsup\o_s(x_ixx_i)\\
 &\le&\o_n(x)+ \limsup \e^{2\t(e_i)}\o_s(e_i^\perp xe_i^\perp)\\
 &\le& \o_n(x)+ \e^2\|\o_s\|\,\|x\|.
 \ee
It thus follows that
 $\d(\o)\le\d(\o_n)$, so $\d(\o)=\d(\o_n)$.
Now it is easy to conclude the validity of the result. Indeed, the
preceding two parts imply
 $$\d(\o)=\inf\{\t(wx)\;:\; x\in\M_+^{-1},\; \De(x)\ge1\}.$$
By a formula on determinants from \cite{arvesonAJM}, the last
infimum is nothing but $\De(w)$. Therefore, the theorem is proved.
\cqd

\begin{rk}
 The proof above shows that  the infimum in
Theorem \ref{szego formula} remains the same if one requires $a$
to be invertible with $a^{-1}\in\A$ (i.e., $a\in\A^{-1}$). Namely,
 \be
 \d(\o)=\inf\{\o(|a|^p)\;:\; a\in\A^{-1}, \;\De(\F(a))\ge1\}
 =\inf\{\o(|a|^p)\;:\; a\in\A^{-1}, \;\De(a)\ge1\}.
 \ee
 \end{rk}

\medskip

\n{\bf Acknowledgements.} We thank David Blecher and Louis
Labuschagne for keeping us informed of their recent works on
noncommutative Hardy spaces.



\begin{thebibliography}{11}

\bibitem{arvesonAJM}
W.~B. Arveson.
\newblock Analyticity in operator algebras.
\newblock {\em Amer. J. Math.}, 89:578--642, 1967.

\bibitem{blouter}
 D.P. Blecher and L.E. Labuschagne.
\newblock Applications of the Fuglede-Kadison determinant: Szeg\"o's
theorem and outers for noncommutative $H^p$.
\newblock Preprint, 2006.

\bibitem{blcharac}
 D.P. Blecher and L.E. Labuschagne.
\newblock Characterizations of noncommutative $H^\8$.
\newblock {\em Integr. Equ. Oper. Theory,} to appear.

\bibitem{blbeurling}
 D.P. Blecher and L.E. Labuschagne.
\newblock A Beurling theorem for noncommutative $L^p$.
\newblock {\em J. Oper. Theory,} to appear.

\bibitem{brown}
 L.G. Brown.
\newblock Lidskii's theorem in the type II case.
\newblock {\it Geometric methods in operator algebras} (Kyoto 1983),
H. Araki and E. Effros (Eds.),  Pitman Res. Notes in Math. Ser.
123, Longmam Sci. Tech. (1986), pages 1--35.

\bibitem{exel}
R. Exel.
\newblock Maximal subdiagonal algebras.
\newblock {\em Amer. J. Math.}, 110:775--782, 1988.

\bibitem{fug-kad}
    B. Fuglede  and R.V. Kadison.
 \newblock Determinant theory in finite factors.
  \newblock {\em Ann. Math.} 55:520--530, 1952.

\bibitem{garn}
    J.B. Garnett.
 \newblock {\em Bounded analytic functions.}
  \newblock Academic Press, 1981.

\bibitem{haagsch}
U. Haagerup and H. Schultz.
\newblock Brown measures of unbounded operators affiliated with
a finite von Neumann algebra.
\newblock Preprint.

\bibitem{kol}
 A. N.~Kolmogoroff.
\newblock   Stationary sequences in Hilbert space.
\newblock {\em Bull. M.G.U.},  Vol. II, 1941.

\bibitem{lab-szego}
L.E. Labuschagne.
\newblock A noncommutative {S}zeg\"o theorem for
subdiagonal subalgebras of von
  {N}eumann algebras.
\newblock {\em Proc. Amer. Math. Soc.}, 133:3643--3646, 2005.


\bibitem{lab-comp}
L.E. Labuschagne.
\newblock Analogues of composition operators on non-commutative $H^p$ spaces.
\newblock {\em J. Operator Theory}, 49:115-141, 2003.


\bibitem{mars-west}
M.~Marsalli and G.~West.
\newblock Noncommutative {$H\sp p$} spaces.
\newblock {\em J. Operator Theory}, 40:339--355, 1998.

\bibitem{px-survey}
G.~Pisier and Q.~Xu.
\newblock Non-commutative {$L\sp p$}-spaces.
\newblock In {\em Handbook of the geometry of Banach spaces, Vol.\ 2}, pages
  1459--1517. North-Holland, Amsterdam, 2003.

\bibitem{saito79}
K-S. Saito.
\newblock A note on invariant subspaces for finite maximal subdiagonal
  algebras.
\newblock {\em Proc. Amer. Math. Soc.}, 77:348--352, 1979.

\bibitem{szego}
 G.~Szeg\"{o}.
\newblock   Beitr\"{a}ge zur Theorie der toeplitzen Formen (Erste
Mitteilung).
\newblock {\em Math. Z.}, 6:167-202, 1920.

\bibitem{tak-I}
M.~Takesaki.
\newblock {\em Theory of operator algebras. {I}}.
\newblock Springer-Verlag, New York, 1979.

\bibitem{ver}
S.~Verblunsky.
\newblock On positive harmoninc functions.
\newblock {\em Proc. London Math. Soc.}, 40:290-320, 1936.

\end{thebibliography}

\end{document}